\PassOptionsToPackage{lowtilde}{url} 
\documentclass[a4paper]{siamart}

\usepackage[utf8]{inputenx}      
\usepackage[T1]{fontenc}         
\usepackage{grffile}             
\usepackage{euscript,graphicx,epsfig,amsfonts,mathtools,epstopdf}
\usepackage[ruled,lined,vlined,linesnumbered,algo2e]{algorithm2e}
\usepackage{siunitx}
\usepackage{listings}
\usepackage{subcaption}
\usepackage{xcolor}
\definecolor{darkgreen}{rgb}{0.0,0.6,0.0}
\usepackage[draft]{changes}

\def\cal{\mathcal}

\def\BA{\textbf{A}}

  \def\CP{{\cal P}}

\def\BAs{\BA{\kern-1.5pt}}

\def\CPs{\CP{\kern-0.8pt}}

\mathcode`\@="8000
{\catcode`\@=\active \gdef@{\mkern1mu}}

\def\mydate{\number\day\ {\ifcase\month \or January\or February\or
              March\or April\or May\or June\or July\or August\or
              September\or October\or November\or December\fi}
\number\year}

\newcommand\restr[2]{{
  \left.\kern-\nulldelimiterspace 
  #1 
  \vphantom{\big|} 
  \right|_{#2} 
  }}

\def\be{\begin{equation}}
\def\ee{\end{equation}}
\def\bea{\begin{eqnarray}}
\def\eea{\end{eqnarray}}

\def\bbmat{\begin{bmatrix}}
\def\ebmat{\end{bmatrix}}

\def\bthm{\begin{theorem}}
\def\ethm{\end{theorem}}
\def\blem{\begin{lemma}}
\def\elem{\end{lemma}}
\def\bprop{\begin{proposition}}
\def\eprop{\end{proposition}}
\def\bcor{\begin{corollary}}
\def\ecor{\end{corollary}}
\def\bdefin{\begin{definition}}
\def\edefin{\end{definition}}
\def\bc{\begin{cases}}
\def\ec{\end{cases}}

\newtheorem{exple}{Example}
\def\bex{\begin{exple}}
\def\eex{\end{exple}}
\def\bass{\begin{assumption}}
\def\eass{\end{assumption}}

\DeclareMathSymbol{\dprod}{\mathbin}{operators}{"3A}

\providecommand{\file}[1]{\texttt{\nolinkurl{#1}}}

\newsiamthm{example}{Example}
\newsiamthm{rem}{Remark}

\renewcommand{\theequation}{\thesection.\arabic{equation}} 
\renewcommand{\thefigure}{\thesection.\arabic{figure}}
\renewcommand{\thetable}{\thesection.\arabic{table}}
\crefname{algocf}{algorithm}{algorithms}
\Crefname{algocf}{Algorithm}{Algorithms}
\numberwithin{theorem}{section}

\newcommand{\TheTitle}{A fast and memory-efficient spectral Galerkin scheme for distributed elliptic optimal control problems} 
\newcommand{\TheAuthors}{L.H. Christiansen, and J.B. J{\O}rgensen}

\headers{A SPECTRAL GALERKIN SCHEME FOR OPTIMAL CONTROL}{\TheAuthors}

\title{{\TheTitle}\thanks{Submitted for review 25/8/2017.}}

\author{
	Lasse H. Christiansen \thanks{Department of Applied Mathematics and Computer Science \& Center for Energy Resources Engineering,
 			Technical University of Denmark (DTU),
 			DK-2800 Kgs. Lyngby, Denmark (\email{lhch@dtu.dk})}
	\and
	John B. J{\O}rgensen \thanks{Department of Applied Mathematics and Computer Science \& Center for Energy Resources Engineering,
 			Technical University of Denmark (DTU),
 			DK-2800 Kgs. Lyngby, Denmark. (\email{jbjo@dtu.dk})}
}

\ifpdf
\hypersetup{
	pdftitle={\TheTitle},
	pdfauthor={\TheAuthors}
}
\fi

\begin{document}

\maketitle

\begin{abstract}
Many scientific and engineering challenges can be formulated as optimization problems which are constrained by partial differential equations (PDEs). These include inverse problems, control problems, and design problems. As a major challenge, the associated optimization procedures are inherently large-scale. To ensure computational tractability, the design of efficient and robust iterative methods becomes imperative. To meet this challenge, this paper introduces a fast and memory-efficient preconditioned iterative scheme for a class of distributed optimal control problems governed by convection-diffusion-reaction (CDR) equations. As an alternative to low-order discretizations and Schur-complement block preconditioners, the scheme combines a high-order spectral Galerkin method with an efficient preconditioner tailored specifically for the CDR application. The preconditioner is matrix-free and can be applied within linear complexity where the proportionality constant is small. Numerical results demonstrate that the preconditioner is ideal in the sense that appropriate Krylov subspace methods converge within a low number of iterations, independently of the problem size and the Tikhonov regularization parameter.
\end{abstract}

\begin{keywords}
  PDE-constrained optimization, Saddle-point problems, Preconditioning, Spectral methods, Optimal control.
\end{keywords}

\begin{AMS}{49J20, 65F08, 65F10, 93C20}
\end{AMS}

\section{Introduction}
PDE-constrained optimal control is ubiquitous in science and engineering. Applications include diverse areas such as geoscience, chemical process industry, aerodynamics, quantum systems, medicine, manufacturing, and financial engineering \cite{leugering2014trends}. However, real-life applications are impeded by computational challenges tied to the large-scale nature of PDE-constrained optimal control. Often, the computational bottleneck amounts to solving large linear saddle-point systems \cite{herzog2010a}. The large dimensions imply that direct linear solvers become inefficient in terms of both CPU-time and storage. Consequently, proper design of preconditioned iterative methods are of crucial importance. 
The literature provides a variety of preconditioned schemes targeting specific applications. These include Poisson control \cite{rees2010a,schoeberl2007a}, Stokes control \cite{rees2011a}, Navier-Stokes control \cite{dolgov2017a,pearson2015a}, and optimal control of reaction-diffusion processes \cite{pearson2013a}. Most of these schemes combine low-order finite element discretizations with Schur-complement block preconditioners. While this constellation has proven effective for a range of model problems, computational efficiency relies heavily on the availability of inexpensive approximations to the Schur-complement. \par As an alternative to existing methodologies, this paper proposes a fast, high-precision and memory-efficient spectral Galerkin (SG) scheme tailored for distributed elliptic optimal control problems governed by convection-diffusion-reaction (CDR) equations. Such problems are of importance to e.g. chemical process engineering, mathematical biology and simulation of fluid flow. Compared to existing practices, the SG scheme contributes in two ways. Firstly, in contrast to the wide use of low-order discretizations, the SG scheme uses a custom-made \textsl{high-order spectral Galerkin method}. This is motivated by the observation that when 1) the optimization problem is constrained only by the PDE or when 2) integral constraints are imposed on the state - and/or control variable, spectral methods converge \textsl{spectrally} fast to the optimal solution, provided only that the data of the problem is smooth \cite{chen2008a}. Consequently, compared to low-order finite difference or finite element methods, the SG scheme requires significantly fewer unknowns to resolve the problem to a given precision. Effectively, this reduces computational complexity by lowering both storage requirements and CPU-time. Secondly, the scheme uses a new preconditioner tailored specifically to CDR optimal control problems that arise from the SG discretization. At its core, the preconditioner uses an efficient direct solver for linear diffusion-reaction type equations with constant coefficients. In this way, it resembles a fast Poisson preconditioner known from traditional boundary value problems \cite{concus1973a}. Unlike most preconditioners that target PDE-constrained optimization, it does not rely on the availability of Schur-complement approximations. Further, the preconditioner is matrix-free and can be applied within linear complexity where the proportionality constant is small.  
\par Previously, efficient spectral Galerkin schemes have been developed in the context of Poisson control \cite{chen2011a,chen2008a}. Here the authors proposed an efficient direct solver that relies on diagonalization \cite{shen1994a}.  As a drawback, the methods are restricted to linear constant coefficient problems. In this perspective, this paper provides a non-trivial extension of \cite{chen2011a, chen2008a} to efficiently handle more complex cases of variable-coefficients. In this way, this research may pave the way for new applications of high-order methods to PDE-constrained optimal control.
To demonstrate the schemes potential, two case studies solve CDR distributed optimal control problems in 2D and 3D, respectively. The numerical results show that spectrally accurate optimal controls can be obtained within a few iterations of a suitable preconditioned Krylov subspace method. Further, the results show that the iteration count is independent of both the problem size and the Tikhonov regularization parameter. This strongly suggest that the preconditioner is ideal for the application. \par
The paper is organized as follows. Section \ref{Sec2} introduces the optimal control problem and derives first-order necessary optimality conditions. Discretization of the optimality conditions is described in Section \ref{Sec3}. Section \ref{Sec4} introduces the SG scheme, including the new preconditioner. Section \ref{Sec5} presents numerical results and conclusions are made in Section \ref{Sec6}.

\section{A class of elliptic optimal control problems} \label{Sec2}
This paper treats the class of distributed elliptic optimal control problems in the form: 
\begin{subequations} \label{OCPStrong1}
\begin{alignat}{5} 
&\min_{y \in Y, \hspace{2pt} u \in U}  \quad  && \frac{1}{2} \Vert y-y_{d} \Vert_{L^2(\Omega)}^2 + \frac{\rho}{2} \Vert u \Vert_{L^2(\Omega)}^2, \hspace{3pt} \rho >0,  \label{ObjStrong} \\
&  \quad  \text{s.t.} \quad &&  - \nabla \cdot \mathcal{F}(y)+ \gamma y = u+f  \quad && \text{in} \quad && \Omega.  \label{PDEStrong1}   \\
&   \quad && y= g \quad && \text{on} \quad && \Gamma. \label{PDEStrong2}
\end{alignat}
\end{subequations}
Here the flux, $\mathcal{F},$ is given by $\mathcal{F}(y)= \alpha \nabla y + \beta y$. 
 \newline \newline  Tracking-type PDE-constrained problems of the type (\ref{OCPStrong1}) trace back to the seminal work of Lions \cite{lions1971optimal}. Such problems seek to control the state variable, $y \in Y$, to a pre-specified set-point, $y_{d} \in Z,$ by manipulating the source of control, $u \in U$.  The objective (\ref{ObjStrong}) expresses the desire to reach the set-point, $y_{d}$, while the PDE-constraints, (\ref{PDEStrong1}) and (\ref{PDEStrong2}), ensure that the optimal solution, $(y^{\ast},u^{\ast}) \in Y \times U,$ satisfies the underlying physical model, where $g$ specifies the state on the boundary. \newline The constant $\rho$ denotes the Tikhonov regularization parameter. PDEs of the type (\ref{PDEStrong1}) model convection-diffusion-reaction processes. The equation typically derives from first principle mass conservation laws and can be written in the equivalent form
 \begin{equation} L(y) \equiv \displaystyle \sum_{i,j=1}^{d} D_{i} (\alpha D_{j} y) + \displaystyle \sum_{i=1}^{d} D_{i}(\beta_{i} y )+ \gamma y=u+f.  \label{ExplicitForm} \end{equation}
 Here $D_{i}$ denotes the partial derivative, $\frac{\partial}{\partial x_{i}}, \hspace{3pt} \alpha(\cdot)$ denotes the diffusion coefficient, $\beta(\cdot)$ is a vector that specifies the direction and magnitude of convection and $\gamma(\cdot)$ describes reaction.  
 Fundamentally, there exists two approaches for the numerical solution of (\ref{OCPStrong1}). The literature distinguishes between the strategies of 1) \textsl{optimize-then-discretize} (OD) and 2) \textsl{discretize-then-optimize} (DO). The OD approach uses techniques of optimization in Banach spaces to derive first-order necessary optimality conditions and solves the associated infinite-dimensional system using an appropriate discretization scheme. In contrast, the DO approach first discretizes (\ref{OCPStrong1}) by replacing the objective (\ref{ObjStrong}) and the PDE-constraints, (\ref{PDEStrong1}-\ref{PDEStrong2}), by finite dimensional counterparts. In this way, the DO strategy allows for the use of conventional finite-dimensional optimization theory \cite{Nocedal}.  
In general, DO and OD do \textsl{not} commute \cite{leykekhman2012a}, i.e., the approaches may lead to different finite-dimensional optimization problems depending on the discretization scheme. 
In general, there is no broad consensus on the preferred method and the choice between DO and OD often depends on the specific application. For a more detailed discussion, see e.g. \cite{hinze2008optimization}.  To solve (\ref{OCPStrong1}), this paper uses the OD approach. This choice allows for construction of numerical examples that can verify convergence of the SG scheme to the true infinite-dimensional optimal state and control, $(y^{\ast},u^{\ast}) \in Y \times U$.
\subsection{OD - Variational formulation and first-order optimality conditions}
Following the OD approach, this section introduces the variational formulation of (\ref{OCPStrong1}), specifies underlying assumptions on the problem data and derives first-order necessary optimality conditions that characterize the optimal solution $(y^{\ast},u^{\ast}) \in Y \times U$.
To this end, let $\Omega \subset \mathbb{R}^{d}, \hspace{3pt} d=1,2,3$ be a bounded rectangular domain with boundary, $\Gamma.$ Without loss of generality, the following assumes homogeneous boundary conditions, i.e., $g:=0$. Inhomogeneous boundary conditions are readily accounted for by appropriately modifying the source, $f.$ Take $Y:=H_{0}^{1}(\Omega)$ to be the state space, let $U:=L^2(\Omega)$ be the set of admissible controls and consider the space of test functions $V:=H^1_{0}(\Omega).$ The variational formulation of the optimal control problem (\ref{OCPStrong1}) then becomes:
\begin{subequations} \label{OCPWeak}
\begin{alignat}{5}
&\min_{y \in Y, \hspace{2pt} u \in U}  \quad  &&  \frac{1}{2} \Vert y-y_{d} \Vert_{L^2(\Omega)}^2 + \frac{\rho}{2} \Vert u \Vert_{L^2(\Omega)}^2, &&\hspace{3pt} \rho >0, \\
&  \quad  \text{s.t.} \quad &&  a(y,v) - b(u,v) =  \langle f,v \rangle_{L^2(\Omega)} \quad && \forall \hspace{3pt} v \hspace{3pt} \in V.
\end{alignat} 
\end{subequations}
Here the bilinear forms $a: Y \times V \rightarrow \mathbb{R}, \hspace{3pt} b: U \times V \rightarrow \mathbb{R}$  are given by
\begin{align}
&a(y,v) = \displaystyle \int_{\Omega} (\mathcal{F}(y) \cdot \nabla v+ \gamma y v) dx \\
& b(u,v) =  \displaystyle \int_{\Omega} uv dx.
\end{align}
To ensure well-posedness of the problem (\ref{OCPWeak}), this paper assumes that $ \alpha , \gamma \in L^{\infty}(\Omega)$ are strictly positive and that $\beta_{i} \in W^{1,\infty}(\Omega), \hspace{3pt} f, y_{d} \in L^2(\Omega).$ Also, this paper restricts attention to diffusion-dominated problems, i.e.,  $|\beta|$ is assumed to be of moderate size. For convection-dominated problems, see e.g. \cite{HH2013}. Lastly, it is assumed that the following inequality holds 
\begin{equation} - \frac{1}{2} \nabla \cdot \beta + \gamma \geq 0, \hspace{3pt} \text{in} \hspace{3pt} \Omega. \end{equation}
Under these assumptions, it is well-known that the optimal control problem (\ref{OCPWeak}) admits a unique solution $(y^{\ast}, u^{\ast}) \in Y \times U$ \cite{lions1971optimal,rotzsch2010optimal}.
Further, optimization theory in Banach spaces shows that the unique optimal solution to (\ref{OCPWeak}) corresponds to a stationary point of the Lagrangian function $\mathcal{L}: H_{0}^{1}( \Omega) \times L^2(\Omega) \times H_{0}^{1}(\Omega) \rightarrow \mathbb{R}$ \cite{rotzsch2010optimal}: 
\begin{equation}  \mathcal{L}(y,u,p):=  \frac{1}{2} \Vert y-y_{d} \Vert_{L^2(\Omega)}^2 + \frac{\rho}{2} \Vert u \Vert_{L^2(\Omega)}^2 + a(y,p) - b(u,p) - \langle f,p \rangle_{L^2(\Omega)}.\end{equation} 
This is equivalent to the requirement that all partial Fr\'echet derivatives of the Lagrangian vanish, i.e., the unique solution $(y^{\ast}, u^{\ast}) \in Y \times U$ to (\ref{OCPWeak}) must satisfy that $\nabla \mathcal{L}(y^{\ast},u^{\ast},p^{\ast}) =0,$ where $p^{\ast} \in H^1_{0}(\Omega)$ denotes the corresponding adjoint state. This leads to the set of necessary and sufficient optimality conditions:
\begin{subequations} \label{OPT}
\begin{alignat}{5} 
&a(y^{\ast},v) + b(u^{\ast},v) =  \langle f,v \rangle_{L^2(\Omega)} \quad && \forall \hspace{3pt} v \hspace{3pt} \in V, \quad && \text{(State equation)}\\
& a(v,p^{\ast}) =\langle y_{d}-y^{\ast},v \rangle_{L^2(\Omega)}, \quad && \forall \hspace{3pt} v \in V, \quad && \text{(Adjoint equation)} \\
& b(w,p^{\ast}) = \rho \langle u^{\ast},w \rangle,  \quad &&  \forall \hspace{3pt} w \in U,  \quad && \text{(Gradient equation)}. \label{OPTGrad}
\end{alignat} 
\end{subequations} 
Note that the gradient equation (\ref{OPTGrad}) leads to the condition 
\begin{equation}  \displaystyle \int_{\Omega} (\rho u^{\ast}-p^{\ast}) w dx =0, \hspace{3pt} \forall \hspace{3pt} w \in U.\end{equation} 
Since $U=L^2(\Omega)$, it follows that $\rho u^{\ast}(x) = p^{\ast}(x), \hspace{3pt} \text{a.e.} \hspace{3pt} x \in \Omega.$ Consequently, the optimality conditions (\ref{OPT}) reduce to the dual system of coupled PDEs
\begin{subequations} \label{dualsys}
\begin{alignat}{5}
&a(y^{\ast},v) - b(u^{\ast},v) =  \langle f,v \rangle_{L^2(\Omega)}  \quad && \forall \hspace{3pt} v \in V, \quad && \text{(State equation)}\\
& \rho a(v,u^{\ast}) +b(y^{\ast},v) =\langle y_{d},v \rangle_{L^2(\Omega)}, \quad && \forall \hspace{3pt} v \in V, \quad && \text{(Adjoint equation)}. 
\end{alignat} 
\end{subequations} 
 \section{OD - Discretization of the optimality system} \label{Sec3}
The following describes how to discretize the optimality system (\ref{dualsys}) using a customized spectral Galerkin method. While the method was first introduced by \cite{shen2007a} in the context of traditional initial- boundary value problems, this paper is the first to recognize its computational benefits to PDE-constrained optimal control. The method constitutes a corner stone of the SG scheme and plays an integral role in the design of the new preconditioner. Further, the spectral Galerkin method converges \textsl{exponentially} fast to the optimal solution, $(y^{\ast}, u^{\ast}) \in Y \times U,$ provided only that $f$ and $z_{d}$ are smooth \cite{chen2008a,zhou2015a}. In turn, spectral accuracy implies that only a few number of expansion modes, $N$, are required to reach a given level of accuracy. Effectively, this lowers dimensionality of the problem. For the class of CDR problems (\ref{OCPStrong1}), the spectral Galerkin method therefore provides an efficient alternative to conventional low-order finite differences and finite element discretizations that currently predominate numerical schemes for PDE-constrained optimization.  \subsection{The spectral Galerkin method} \label{SGmethod}
For ease of treatment, the presentation covers only the one-dimensional case on the reference domain $\Omega=[-1,1]$. Multi-dimensional problems can be treated accordingly by forming tensor products of the one-dimensional components. 
To discretize the optimality system (\ref{dualsys}), let $\mathbb{P}_{N}$ denote the set of all polynomials of degree less than or equal to $N$ and consider the subspace
\begin{equation} V_{N}= \lbrace v \in \mathbb{P}_{N} :  v( \pm 1)=0 \rbrace.\end{equation} 
Note that $V_{N}$ is a subspace of both $Y=V:=H^{1}_{0}(\Omega)$ and $U:=L^2(\Omega).$ 
The spectral Galerkin method approximates the true solution $(y^{\ast},u^{\ast}) \in Y \times U$ by the truncated series expansions
\begin{align}
y_{N}(x) &= \displaystyle \sum_{k=0}^{N-2} \widehat{y}_{k} \psi_{k}(x), \hspace{3pt}  \bold{\widehat{y}} =\lbrace  \widehat{y}_{k}\rbrace_{k=0}^{N-2} , \label{ExpState} \\
u_{N}(x) &= \displaystyle \sum_{k=0}^{N-2} \widehat{u}_{k} \psi_{k}(x), \hspace{3pt}  \bold{\widehat{u}} =\lbrace  \widehat{u}_{k}\rbrace_{k=0}^{N-2}. \label{ExpControl}
\end{align}
Here $\lbrace \phi_{k} \rbrace_{k=0}^{N-2}$ is a suitable basis for $V_{N}$  and $\bold{\widehat{y}}, \hspace{3pt} \bold{\widehat{u}}$ denote the unknown expansion coefficients. 
The choice of basis represents one of the primary differences between spectral - and finite element methods. While finite element discretizations employ local bases with finite regularity, spectral methods rely on global and smooth basis functions such as Chebyshev, Legendre, Laguerre and Hermite polynomials \cite{canuto2007spectral}. The efficiency of spectral discretizations rely crucially on the choice of basis. In this regard, this paper uses a set of Fourier-like basis functions that was first introduced by \cite{shen2007a}: 
\begin{equation} \psi_{k}(x) = \displaystyle \sum_{j=0}^{N-2} q_{jk} \phi_{j}(x) , \hspace{3pt} 0 \leq k \leq N-2. \label{Fourier-like}\end{equation}
Here $Q=(q_{jk})$ is the orthogonal matrix whose columns are the eigenvectors of the matrix, $\widehat{M}=(\widehat{m}_{ij})$, defined by
\begin{equation} \label{MassMatrix} \widehat{m}_{jk} = \begin{cases} c_{k} c_{j} \left(\frac{2}{2j+1} + \frac{2}{2j+5} \right), & k=j\\
 -c_{k}c_{j} \frac{2}{2k+5}, & k=j+2 
\end{cases}, \hspace{5pt} c_{j} = \frac{1}{\sqrt{4j+6}}, \end{equation}
while the auxiliary functions, $\phi_{j},$ are given by
\begin{equation} \phi_{j}(x) =  c_{j} \left( L_{j}(x) - L_{j+2}(x) \right), \hspace{3pt} 0 \leq j \leq N-2. \label{Legendre} \end{equation}
$\lbrace L_{j}(\cdot) \rbrace_{k=0}^{N}$ denote the orthogonal Legendre polynomials. The benefits of the Fourier-like basis to PDE-constrained optimization will be elaborated in Section \ref{FastDirect}.
To determine the associated unknown expansion coefficients, $\bold{\widehat{y}}, \hspace{3pt} \bold{\widehat{u}}$, the spectral Galerkin method requires that the discrete optimality conditions hold
\begin{subequations} \label{dualsysdis}
\begin{alignat}{5}
&a(y_{N},v_{N}) - b(u_{N},v_{N}) =  \langle I_{N} f,v_{N} \rangle_{L^2(\Omega)}  \quad && \forall \hspace{3pt} v_{N} \hspace{3pt} \in V_{N}, \\
& \rho a(w_{N},u_{N}) + b(y_{N},w_{N}) =\langle I_{N}y_{d},w_{N} \rangle_{L^2(\Omega)}, \quad && \forall \hspace{3pt} w_{N} \in U_{N}.
\end{alignat} 
\end{subequations}
Here  $I_{N}$ denotes the interpolation operator over the set of Legendre-Gauss-Lobatto (LGL) nodes, $\lbrace x_{j} \rbrace_{j=0}^{N}$. For details on LGL nodes and numerical quadrature, see e.g. \cite{kopriva2009implementing}.
Since $Q$ has full rank, the set of functions $\lbrace \psi_{k} \rbrace_{k=0}^{N-2}$ form a new basis for $V_{N}.$ Hence, the conditions (\ref{dualsysdis}) are equivalent to the coupled system of equations
\begin{subequations} \label{linearsysdis}
\begin{alignat}{5}
&a(y_{N},\psi_{k}) - b(u_{N},\psi_{k}) =  \langle I_{N} f, \psi_{k} \rangle_{L^2(\Omega)}  \quad && 0 \leq k \leq N-2,\\
& \rho a(\psi_{k},u_{N}) + b(y_{N}, \psi_{k}) =\langle I_{N}y_{d}, \psi_{k} \rangle_{L^2(\Omega)}, \quad && 0 \leq k \leq N-2.
\end{alignat} 
\end{subequations}
Using the notation, 
\begin{alignat}{5}
& y_{d_{k}}= \langle I_{N} y_{d}, \psi_{k} \rangle, \quad && Y= (y_{d_{0}},y_{d_{1}}, \hdots, y_{d_{N-2}})^{T}, \\
& f_{k}= \langle I_{N} f, \psi_{k} \rangle, \quad && F = (f_{0},f_{1}, \hdots, f_{N-2})^{T},\\
& b_{ij} = a(\psi_{j}, \psi_{i}), \quad && B= (b_{ij})_{i,j=0..N-2}, \label{Bmat} \\
& m_{ij} = b(\psi_{j}, \psi_{i}), \quad && M= (m_{ij})_{i,j=0..N-2}, \label{Mmat}
\end{alignat}
the discrete optimality system (\ref{linearsysdis}) can be written in matrix form
\begin{equation} \label{SaddlePointSys}
\underbrace{ \left[ \begin{array}{cc}
    M& \rho B^{T} \\ 
    B & -M \\ 
  \end{array} \right]}_{A} \underbrace{\left[   \begin{array}{c}
    \widehat{y} \\ 
    \widehat{u} \\
  \end{array} \right]}_{x} = \underbrace{ \left[   \begin{array}{c}
    Y \\ 
    F \\
  \end{array} \right]}_{b}.
\end{equation}
For a more detailed introduction to spectral discretization methods and their applications, see e.g. the monographs \cite{canuto2007spectral,hesthaven2007spectral,shen2011spectral} and the list of references therein.  
\section{The SG scheme - An efficient preconditioned iterative solver} \label{Sec4}
The discrete optimality conditions (\ref{SaddlePointSys}) constitute a linear system in saddle-point form \cite{benzi2005a}. Here each of the blocks $(1,2)$ and $(2,1)$ represent a discretized PDE. Consequently, the system is inherently large-scale. For constant-coefficient problems, there exists efficient direct solvers \cite{chen2011a}. However, for variable-coefficient problems, the global nature of the spectral basis functions implies that the matrix $B$ becomes dense. Hence, inversion of the matrix requires $O(N^{3d})$ operations while formation and storage takes additional $O(N^{2d})$ operations. As a result, direct solution strategies become computationally demanding even for a moderate number of unknowns. Despite spectral accuracy, this prevents practical applications in three dimensions. To meet this challenge, the following introduces the SG scheme that is tailored specifically for the efficient solution of systems (\ref{SaddlePointSys}) that arise from the spectral Galerkin discretization. As an addition to \cite{chen2011a}, the SG scheme provides an efficient way to handle variable-coefficients. At its foundation, the scheme uses non-symmetric Krylov subspace (KSP) methods, such as GMRES \cite{saad1986a} and BiCG \cite{Fletcher1976}. Fundamentally, KSP methods solve the linear system (\ref{SaddlePointSys}) iteratively by approximating the true solution, $x,$ using a sequence of iterates, $x_{k},$ that are extracted from the $k$-dimensional Krylov subspaces:
\begin{equation} \mathcal{K}(A,r_{0}) := \text{span} \hspace{3pt} \lbrace r_{0}, Ar_{0}, A^2 r_{0},...,A^{k-1} r_{0} \rbrace.\end{equation} 
Here $r_{0} :=b-Ax_{0}$ denotes the initial residual vector.  
For large-scale problems, the efficiency of KSP methods relies on two key components:
 \begin{itemize}
 \item[(i)] Matrix-free and efficient evaluation of the matrix-vector products, $Ax_{k},$ to avoid forming and storing the matrix, $A$.  
 \item[(ii)] An efficient preconditioner, $P,$ to lower the iteration count of the given KSP method.
 \end{itemize} 
In essence, the SG scheme comprises a new strategy of how to realize the requirements, $ \text{(i)}$ and $\text{(ii)},$ in the context of the optimality system (\ref{SaddlePointSys}). This naturally divides the presentation of the scheme into the subsections \ref{MatrixfreeSec} and \ref{PreSec}.  
  \subsection{Matrix-free evaluation of KSP matrix-vector products} \label{MatrixfreeSec}
It is well-known that  the complexity of matrix-vector products involving $A$ can be reduced to $\mathcal{O}(N^{d+1})$ operations using sum-factorization methods that exploit the tensor-product structure of the block matrices \cite{canuto2007spectral}. The following describes how computational costs can be further reduced using transform methods. For simplicity, the presentation focuses on the one-dimensional case. Multi-dimensional transforms can be performed through a sequence of one-dimensional transforms. Hence, the computational bottleneck is tied to $(N-1) \times (N-1)$ matrix-matrix products. This implies that the computational complexity amounts to a small multiple of $N^{d+1}.$  Combined with the high accuracy of the SG method, this makes the approach competitive, even for large scale problems. Also, since matrix-matrix products constitute the computational bottleneck, the scheme is highly amenable to parallelization. 
 \subsubsection{KSP matrix-vector products by transform methods}
Given an arbitrary vector, $[\widehat{x}_1,  \widehat{x}_2]^{T} \in \mathbb{R}^{2 (N-1)},$ consider the matrix-vector product 
\begin{equation} \label{BlockPreconditioner2}
\underbrace{ \left[ \begin{array}{cc}
  M &  \rho B^{T}\\ 
 B & -M \\ 
 \end{array} \right]}_{A} \hspace{2pt} \underbrace{\left[   \begin{array}{c}
 \widehat{x}_1 \\ 
    \widehat{x}_2 \\
  \end{array} \right]}_{x_{k}}.
 \end{equation}
Section \ref{FastDirect} shows that the use of the Fourier-like basis (\ref{Fourier-like}) implies that $M$ becomes a diagonal matrix.  Hence, the matrix-vector products $M \widehat{x}_{i}, \hspace{2pt} i=1,2$ reduce to Euclidean scalar products. The following therefore focuses on the matrix-vector products that involve the full matrices, $B$ and $B^{T}$. To this end, define 
\begin{equation} u_{N,i} := \displaystyle \sum_{\ell=0}^{N-2} \widehat{x}_{i,\ell} \psi_{\ell} \in V_{N}, \hspace{3pt} i=1,2,\end{equation}
where $\widehat{x}_{i}=\lbrace \widehat{x}_{i,\ell} \rbrace_{\ell=0}^{N-2}.$
Since  $B= (b_{ij})_{i,j=0..N-2},$ where $b_{ij} = a(\psi_{j}, \psi_{i})$ it follows that 
\begin{equation}
 \label{MatrixFreeEval} (B \widehat{x}_{i}) _{j} = a(u_{N,i}, \psi_{j}), \hspace{3pt} j=0,..., N-2. \end{equation}
 Hence, the matrix-vector products that involve $B$ can be computed matrix-free by evaluation of the inner products (\ref{MatrixFreeEval}) that take the explicit form
 \begin{equation} \underbrace{\displaystyle \int_{\Omega}  \alpha \nabla u_{N,i}  \cdot \nabla  \psi_{j}  dx}_{T_{1}} +  \underbrace{ \displaystyle \int_{\Omega}  (\beta u_{N,i})  \cdot \nabla  \psi_{j} dx}_{T_{2}} +  \underbrace{\displaystyle \int_{\Omega}  \gamma u_{N,i}   \psi_{j}  dx}_{T_{3}}, \hspace{3pt} j=0,..., N-2. \end{equation} 
 Similarly, the matrix-vector products, $(B^{T} \widehat{x}_{i}) _{j},$ can be computed matrix-free by evaluation of the inner products $a(\psi_{j},u_{N,i}).$
In both cases, the evaluation can be carried out efficiently using transform methods. For simplicity, consider evaluation of the term $T_{3}:$
 \begin{itemize}
 \item [(i)] Evaluate $u_{N,i}$ at the LGL nodes, $\lbrace x_{j} \rbrace_{j=0}^{N},$ using the forward Legendre transform
      \begin{equation} \label{LTransform} u_{N,i}(x_{j}) = \displaystyle \sum_{k=0}^{N-2} \widehat{x}_{i,k} \psi_{k}(x_{j}) , \hspace{3pt} j=0,1,...,N. \end{equation}
 \item [(ii)]  Compute the coefficients, $ \lbrace \widehat{w} \rbrace_{k=0}^{N-2},$ by the backward Legendre transform 
    \begin{equation} I_{N}(\gamma u_{N,i}) (x_{j}) = \displaystyle \sum_{k=0}^{N-2} \widehat{w}_{k} \psi_{k}(x_{j}), \hspace{3pt} j=1,2,...,N-1. \end{equation}
  \item[(iii)] Evaluate the inner products, $T_{3},$ by computing
  \begin{equation} \label{LastStep} \langle  \gamma u_{N,i}, \psi_{j} \rangle = \displaystyle \sum_{k=0}^{N-2} \widehat{w}_{k} \langle \psi_{k}, \psi_{j} \rangle, \hspace{3pt} j=0,1,...,N-2. \end{equation}
  \end{itemize} 
 The terms $T_1$ and $T_2$ can be treated similarly. Each of the transform steps, $\text{(i)}-\text{(ii)}$, require $O(N^{d+1})$ operations while step $\text{(iii)}$ can be carried out in $O(N^{d})$ operations. In all cases, the proportionality constants are small. Hence, the evaluation scheme is highly efficient for moderate size problems and remains reasonable efficient for most large-scale problems. 
\subsection{The SG preconditioner} \label{PreSec}
The role of a preconditioner, $P$, is to modify the original system (\ref{SaddlePointSys}) to make it amenable to iterative solution strategies:
\begin{equation}  P^{-1} Ax= P^{-1} b.\end{equation} 
To effectively lower the iteration count, $P^{-1}$ must provide a good approximation to $A^{-1}$.  However, to be computationally efficient, the preconditioner must be matrix-free and it should be inexpensive to compute the action of $P^{-1}.$ Good preconditioners manage to properly trade off these opposing objectives. In the context of PDE-constrained optimal control, preconditioned iterative methods have gained considerable attention \cite{battermann1998a,battermann2001a,herzog2010b}. In recent years, most contributions have focused on Schur-complement block preconditioners \cite{pearson2015a,pearson2013a,rees2010a,rees2011a}: 
\begin{equation} \label{BlockPreconditionerSchur}
\mathcal{S} :=  \left[ \begin{array}{cc}
   M & \cdot  \\ 
  B & \widehat{S} \\ 
  \end{array} \right].
  \end{equation}
Here $\widehat{S}$ represents an approximation to the negative Schur-complement $S:=-M+\rho B^{T} M^{-1}B.$ Preconditioners of the type (\ref{BlockPreconditionerSchur}) have proven effective in practice and their spectral properties are well understood for a broad class of model problems. However, to be computationally efficient, it is crucial to have approximations to the Schur-complement that are inexpensive to invert. For low-order finite element discretizations, it has become common practice to approximate $S^{-1}$ using a fixed number of standard multi-grid cycles \cite{herzog2010b,rees2010a}. For such schemes, the multi-grid approach leads to preconditioners with linear complexity in the number of variables. However, for spectral Galerkin methods, multi-grid approaches are far less developed. To promote the benefits of spectral schemes, this motives alternatives that do not rely on multi-grid Schur-complement block preconditioners. To this end, this paper proposes preconditioners in the form
\begin{equation} \label{BlockPreconditioner}
 P_{\rho} :=  \left[ \begin{array}{cc}
   M &  \rho \widehat{B} \\ 
  \widehat{B}  & -M \\ 
  \end{array} \right]. 
   \end{equation}
The idea is to replace the dense block, $B,$ of the original system (\ref{SaddlePointSys}), by the \textsl{symmetric} approximation, $\widehat{B}$, that arises from a spectral discretization of the associated \textsl{constant-coefficient} problem:
\begin{subequations} \label{OCPStrongAux}
\begin{alignat}{5} 
&\min_{y \in Y, \hspace{2pt} u \in U}  \quad  && \frac{1}{2} \Vert y-y_{d} \Vert_{L^2(\Omega)}^2 + \frac{\rho}{2} \Vert u \Vert_{L^2(\Omega)}^2, \hspace{3pt} \rho >0,  \label{ObjStrongAux} \\
&  \quad  \text{s.t.} \quad &&  - \bar{\alpha} \Delta y + \bar{r} y=u+f.   \quad && \text{in} \quad && \Omega,  \label{PDEStrong1Aux}   \\
&   \quad && y= g \quad && \text{on} \quad && \Gamma \label{PDEStrong2Aux},
\end{alignat}
\end{subequations}
where $\bar{\alpha}, \hspace{3pt} \bar{r}$ are appropriately chosen constants. 
In this way, the preconditioners (\ref{BlockPreconditioner}) are constructed from the individual blocks of the optimality system associated with the constant-coefficient problem (\ref{OCPStrongAux}). As such, the preconditioners resemble fast Poisson solvers known from traditional BVPs \cite{concus1973a,shen1996efficient}. 
\begin{rem}
Note that the constant-coefficient discretized PDE, $\widehat{B},$ can also be used to approximate the Schur-complement associated with the discrete optimality system (\ref{SaddlePointSys}), i.e., $\widehat{S}:=-M+\rho \widehat{B} M^{-1}\widehat{B}.$ The constant-coefficient approach can therefore also be used to construct efficient Schur-complement block preconditioners (\ref{BlockPreconditionerSchur}) for the optimal control problem (\ref{OCPWeak}).
\end{rem}
\subsubsection{A fast direct solver for constant-coefficient problems} \label{FastDirect}
To apply the preconditioner (\ref{BlockPreconditioner}), each iteration of the KSP method requires solution of the discrete optimality system associated with (\ref{OCPStrongAux}) \begin{equation} \label{SaddlePointSysConst}
 \underbrace{\left[ \begin{array}{cc}
    M& \rho \widehat{B} \\ 
    \widehat{B} & -M \\ 
  \end{array} \right]}_{P_{\rho}} \underbrace{\left[   \begin{array}{c}
    \widehat{y}^{k} \\ 
    \widehat{u}^{k} \\
  \end{array} \right]}_{z_{k}} = \underbrace{\left[   \begin{array}{c}
    \widehat{Y}^{k} \\ 
    \widehat{F}^{k} \\
  \end{array} \right]}_{Ax_{k}}.
\end{equation}
Here $x_{k}$ is the current iterate and $z_{k}:=P_{\rho}^{-1}Ax_{k}$ is the associated intermediate optimal solution. To this end, the preconditioner uses a fast matrix-free \textsl{direct} solver tailored for the constant-coefficient problem (\ref{OCPStrongAux}). In particular, the direct method provides an efficient way to calculate the action of $P_{\rho}^{-1}$ by solving (\ref{SaddlePointSysConst}) in a matrix-free manner. 
 To derive the direct method, consider the Fourier-like basis that was introduced in Section \ref{SGmethod} : 
\begin{equation} \psi_{k}(x) = \displaystyle \sum_{j=0}^{N-2} q_{jk} \phi_{j}(x) , \hspace{3pt} 0 \leq k \leq N-2. \label{Fourier-like2}\end{equation}
This basis share similarities with the Fourier basis for periodic problems. In particular, due to the orthogonality of $Q$, the functions $\lbrace \psi_{k} \rbrace_{k=0}^{N-2}$ satisfy the relations
\begin{subequations}
\begin{align}
\langle \psi_{i}, \psi_{j} \rangle &=  \displaystyle \sum_{n,l=0}^{N-2} q_{ni} q_{lj} \langle \phi_{n}, \phi_{l} \rangle  = (Q^T \widehat{M} Q)_{ji}= \lambda_{j} \delta_{j,i} \label{Relation1} \\
-\langle \psi_{i}^{ \prime \prime}, \psi_{j} \rangle &=  \displaystyle \sum_{n,l=0}^{N-2} q_{ni} q_{lj} \langle \phi_{n}^{\prime \prime}, \phi_{l} \rangle =  \displaystyle \sum_{n,l=0}^{N-2} q_{ni} \delta_{n,l} q_{lj} =  \delta_{j,i}. \label{Relation2}
\end{align}
\end{subequations}
Here $\lbrace \lambda_{k} \rbrace_{k=0}^{N-2}$ are the eigenvalues associated with the eigenvectors that form $Q.$
Combined with the definition of the stiffness matrix, $\widehat{B},$ (\ref{Bmat}), the relations (\ref{Relation1}) and (\ref{Relation2}) imply that 
\begin{equation} \widehat{b}_{ij}= -\overline{\alpha} \langle \psi^{\prime \prime}_{j}, \psi_{i} \rangle + \overline{\gamma} \langle \psi_{j}, \psi_{i} \rangle = \overline{\alpha} \delta_{i,j}+\overline{\gamma} \lambda_i \delta_{i,j}.\end{equation}
Similarly, by definition of the mass matrix (\ref{Mmat}), it follows that $m_{ij}=\lambda_{i} \delta_{ij}.$
Consequently, the Fourier-like basis ensures that the one-dimensional discretized PDE, $\widehat{B},$ and the mass matrix, $M$, become diagonal matrices. In turn, this leads to efficient direct inversion of (\ref{BlockPreconditioner2}). Further, by forming tensor products of the one-dimensional components, efficient solution procedures carry over to the multi-dimensional cases. Consider first the two-dimensional case. Here the mass- and stiffness matrices can be written in the form
\begin{equation}   M_{2D} =   \left[ \begin{array}{ccc}
   \lambda_0 M & \cdot & \cdot \\ 
    \cdot & \ddots & \cdot \\ 
    \cdot & \cdot &  \lambda_{N-2}M \\ 
  \end{array} \right], \hspace{5pt} \widehat{B}_{2D} =  \left[ \begin{array}{ccc}
  \Sigma_{0} & \cdot & \cdot \\ 
    \cdot & \ddots & \cdot \\ 
    \cdot & \cdot & \Sigma_{N-2} \\ 
  \end{array} \right],
\end{equation} 
where $\Sigma_{n}:=\overline{\alpha} (M + \lambda_{n} I_{N-1})+ \overline{\gamma} \lambda_{n} M.$ Hence, using the notation, $$\widehat{y}_{n}^{k}= \lbrace \widehat{y}_{nm}^{k} \rbrace_{m=0}^{N-2}, \hspace{3pt} \widehat{u}_{n}^{k}= \lbrace \widehat{u}_{nm}^{k} \rbrace_{m=0}^{N-2}, \hspace{3pt} \widehat{Y}_{n}^{k}= \lbrace \widehat{Y}_{nm}^{k} \rbrace_{m=0}^{N-2}, \hspace{3pt} \widehat{F}_{n}^{k}= \lbrace \widehat{F}_{nm}^{k} \rbrace_{m=0}^{N-2},$$ it follows that the two-dimensional optimality system (\ref{SaddlePointSysConst}) can be written as $N-1$ independent linear systems
\begin{equation} \label{Sys12} \left[ \begin{array}{cc}
   \lambda_{n} M & \rho \Sigma_{n} \\ 
 \Sigma_{n} &  -\lambda_{n} M  \\ 
  \end{array}  \right]  
\left[ \begin{array}{c}
    \widehat{y}_{n}^{k} \\ 
    \widehat{u}_{n}^{k} \\ 
  \end{array} \right] =  \left[  \begin{array}{c}
  \widehat{Y}_{n}^{k} \\ 
   \widehat{F}_{n}^{k} \\ 
  \end{array} \right],  \hspace{3pt} 0 \leq n \leq N-2. \end{equation} 
Further, since $M$ is a diagonal matrix, the system (\ref{Sys12}) reduces to $(N-1)^{2}$ independent $2 \times 2$ linear systems 
\begin{equation} \label{Sys2d}  \left[ \begin{array}{cc}
  \lambda_{n} \lambda_{m}    & \rho \sigma_{nm} \\ 
   \sigma_{nm} &  - \lambda_{n} \lambda_{m}  \\ 
  \end{array}  \right]  
\left[ \begin{array}{c}
    \widehat{y}_{nm}^{k} \\ 
    \widehat{u}_{nm}^{k} \\ 
  \end{array} \right] =  \left[ \begin{array}{c}
  \widehat{Y}_{nm}^{k} \\ 
   \widehat{F}_{nm}^{k} \\ 
  \end{array} \right],  \hspace{3pt} 0 \leq n,m \leq N-2, \end{equation} 
  where $\sigma_{nm}:=\overline{\alpha}(\lambda_{n}+\lambda_{m}) + \overline{\gamma} \lambda_{n} \lambda_{m}.$
By similar arguments, the three-dimensional discrete optimality system (\ref{SaddlePointSysConst}) can be reduced to $(N-1)^{3}$ independent $2 \times 2$ linear systems 
\begin{equation} \label{Sys3d}  \left[ \begin{array}{cc}
  \lambda_{n} \lambda_{m} \lambda_{k}  &  \rho \sigma_{nmk}  \\ 
    \sigma_{nmk} &  - \lambda_{n} \lambda_{m} \lambda_{k} \\ 
  \end{array}  \right]    
\left[ \begin{array}{c}
    \widehat{y}_{nmk}^{k} \\ 
    \widehat{u}_{nmk}^{k} \\ 
  \end{array} \right] =   \left[  \begin{array}{c}
\widehat{Y}_{nmk}^{k} \\ 
   \widehat{F}_{nmk}^{k} \\ 
  \end{array} \right], \hspace{3pt} 0 \leq n,m,k \leq N-2, \end{equation}
  where $\sigma_{nmk}:= \overline{\alpha} (\lambda_{n} \lambda_{m}+\lambda_{m} \lambda_{k} + \lambda_{n} \lambda_{k}) + \overline{\gamma} \lambda_{n} \lambda_{m} \lambda_{k}.$ \newline \newline
Using the formulation (\ref{Sys12}), it follows that the solution of two and three dimensional optimality systems (\ref{SaddlePointSysConst}) reduces to solution of $N-1$ and $(N-1)^2$ one dimensional problems, respectively. In particular, the action of $P_{\rho}^{-1}$ can be computed within linear complexity, i.e., in $O(N^d)$ operations. Note also that the $2 \times 2$ systems, (\ref{Sys2d}) and (\ref{Sys3d}), are readily inverted analytically. Hence, no direct solver is needed. As a consequence, the solution procedure is essentially matrix-free. In fact, the solution procedure relies only on the $N-1$ eigenvalues of the matrix $\widehat{M},$ defined by (\ref{MassMatrix}). Due to the penta-diagonal structure of $\widehat{M}$, these can be pre-computed and stored efficiently as part of an \textsl{offline} pre-processing stage at the cost of $O(N^2)$ operations. Finally,  the independent nature of the systems makes the solution procedure amenable to parallel implementation.
 \begin{rem}The direct solver presented by this paper shares similarities with the method in \cite{chen2011a}. In particular, both solvers are inspired by the work of Shen \cite{shen1994a} and rely on diagonalization techniques. However, the method that this paper proposes ensures that the action of $P_{\rho}^{-1}$ can be computed matrix-free and within linear complexity, i.e., in $O(N^{d})$ operations. In contrast, use of the direct method in \cite{chen2011a} would require $O(N^{d+1})$ operations to invert the corresponding preconditioner, $P_{\rho}.$ Since KSP methods require the action of $P_{\rho}^{-1}$ at every iteration, this distinction in complexity is important to maintain efficiency for large-scale problems in three dimensions. Note also, that besides serving as a preconditioner, the direct methods charectarized by (\ref{Sys2d}) and (\ref{Sys3d}) may also be used as efficient solvers for 2D and 3D constant-coefficient optimal control problems (\ref{OCPStrongAux}) in their own right.
\end{rem}
\section{Numerical results}  \label{Sec5}
To demonstrate the potential of the SG scheme, two case studies solve prototype CDR optimal control problems (\ref{OCPWeak}) in 2D and 3D, respectively. To support computational efficiency, case study I investigates the SG scheme in terms of convergence and CPU-times for 1) a growing number of unknowns and 2) decreasing error tolerances. The results are compared to 1) \texttt{MATLABs} state-of-the-art direct solver (DS) and 2) a standard second-order finite difference (FD) scheme with a Schur-complement preconditioner (\ref{BlockPreconditionerSchur}), where $M:=I$ is the identity matrix and $B$ is the FD discretized PDE. The FD scheme is included to demonstrate the computational benefits of high-order methods. In particular, it is not intended to represent state-of-the-art. To demonstrate that the preconditioner, $P_{\rho},$ 
is effective and robust, case study II investigates invariance to 1) the problem dimensionality, 2) the Tikhonov regularization parameter, $\rho,$ and to 3) different variable coefficients. The studies consider two sets of coefficients that are chosen as first- and second-order polynomials, respectively. This is motivated by the frequent use of low-order polynomial approximations to non-linear coefficients. Following \cite{concus1973a}, the constant-coefficients, $\overline{\alpha}$ and $\overline{\gamma}$, that are used in the preconditioner, $P_{\rho},$ are defined according to the heuristic
\begin{equation} \overline{\alpha} = \frac{1}{2} \left( \displaystyle \max_{ \bold{x} \in \Omega} \alpha(\bold{x}) + \min_{ \bold{x} \in \Omega} \alpha(\bold{x})\right), \hspace{3pt} \overline{\gamma} = \frac{1}{2} \left( \displaystyle \max_{ \bold{x} \in \Omega} \gamma(\bold{x}) + \min_{ \bold{x} \in \Omega} \gamma(\bold{x})\right). \end{equation}
All computations are performed in \texttt{MATLAB} (2015b) on a 2.9 GHz Intel processor with 16 GB RAM. To solve (\ref{SaddlePointSys}) iteratively, the SG scheme uses \texttt{MATLAB}s implementation of GMRES with a tolerance of $\varepsilon=10^{-10}$ and a maximum of 100 iterations. Due to its low order nature, the FD scheme uses a tolerance of $\varepsilon=10^{-4}.$ All timings are listed in seconds and $N$ refers to the number of modes in each principal direction. 

\subsection{Case study I - Efficiency of the SG scheme} 
\begin{figure}
\centering
\subcaptionbox{Convergence rates in 2D.}{\includegraphics[width=0.45\textwidth]{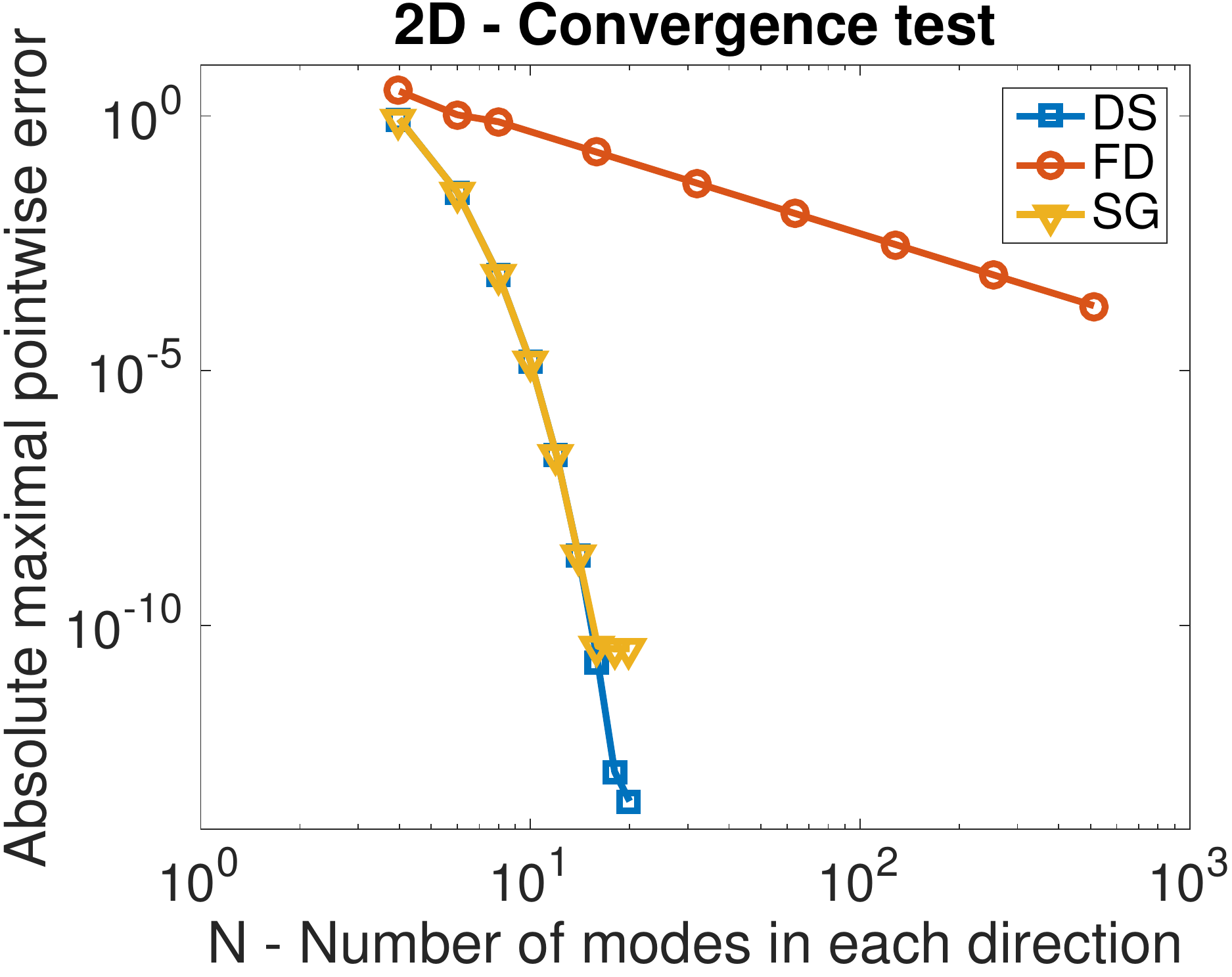}}%
\hfill
\subcaptionbox{Convergence rates in 3D.}{\includegraphics[width=0.45\textwidth]{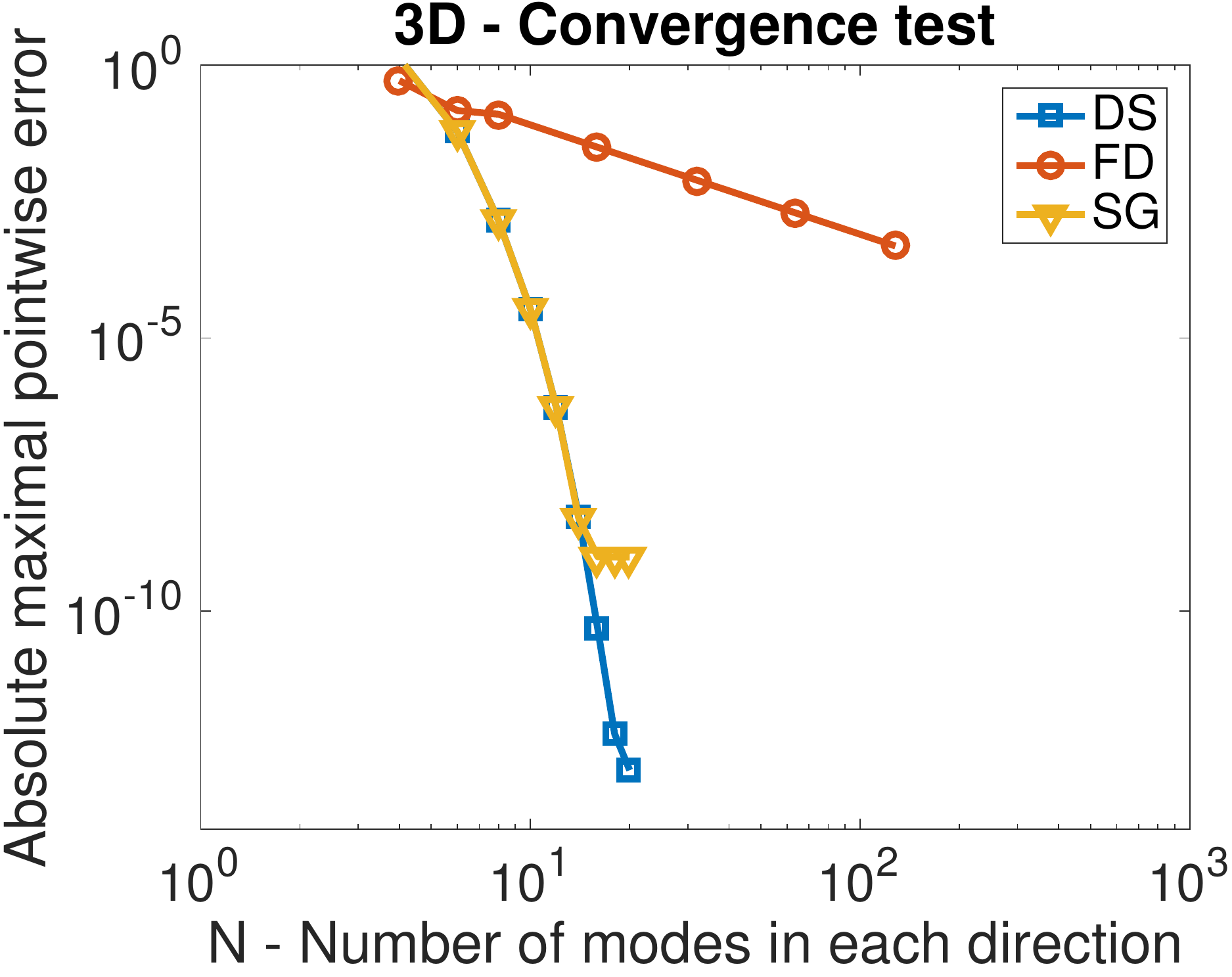}}%
\caption{Convergence rates in 2D and 3D.}
 \label{fig:Convergence}
\end{figure}

\begin{figure}
\centering
\subcaptionbox{2D CPU-times required to reach a given tolerance, $\epsilon$.}{\includegraphics[width=0.45\textwidth]{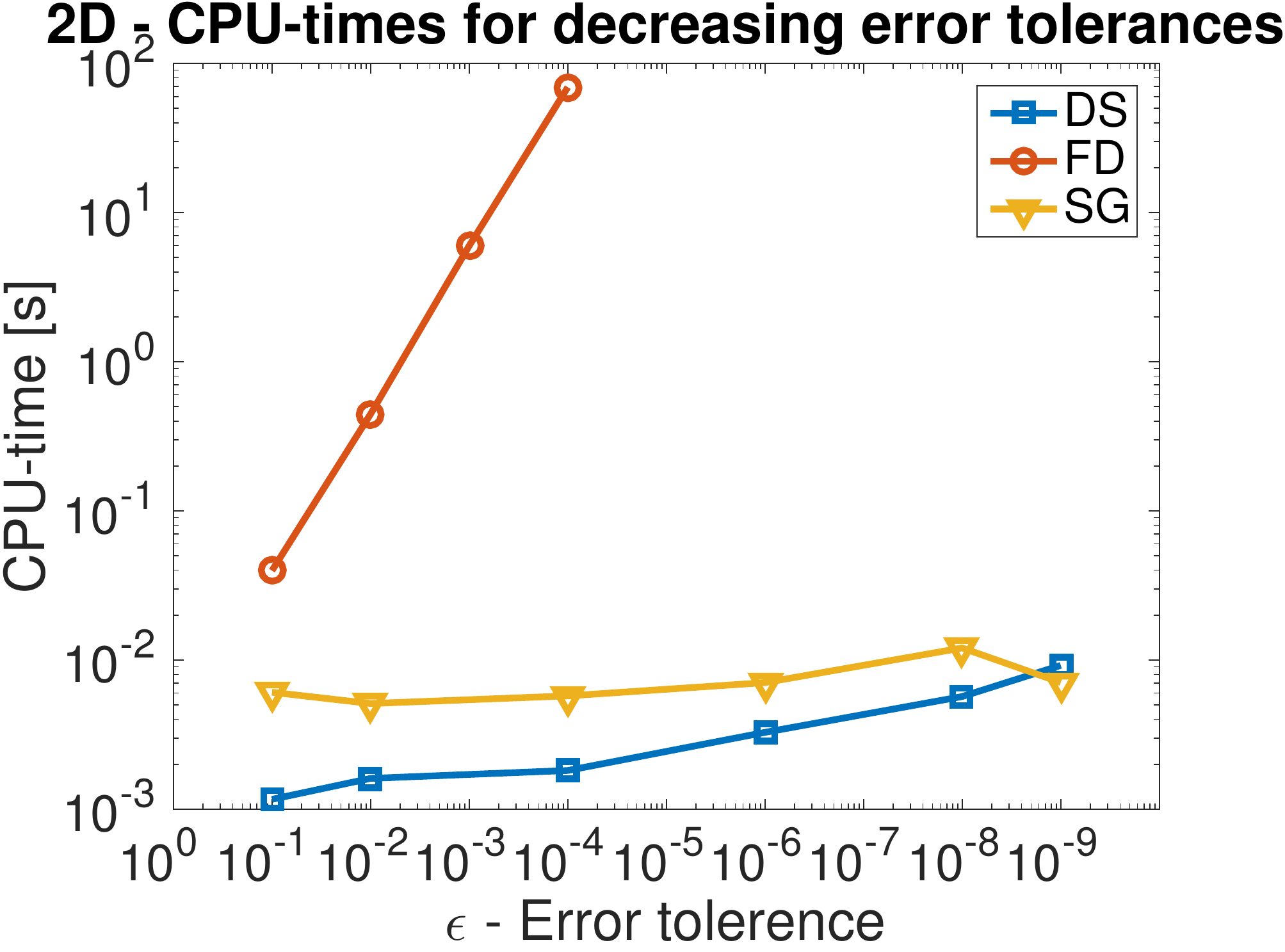}}%
\hfill
\subcaptionbox{3D CPU-times required to reach a given tolerance, $\epsilon$.}{\includegraphics[width=0.45\textwidth]{cpuerror2d.pdf}}%
\caption{CPU-times required to reach a given tolerance, $\epsilon$.}
 \label{fig:CPUtime1}
\end{figure}

\begin{figure}[tb!]
\centering
\subcaptionbox{2D CPU-times for increasing problem sizes.}{\includegraphics[width=0.45\textwidth]{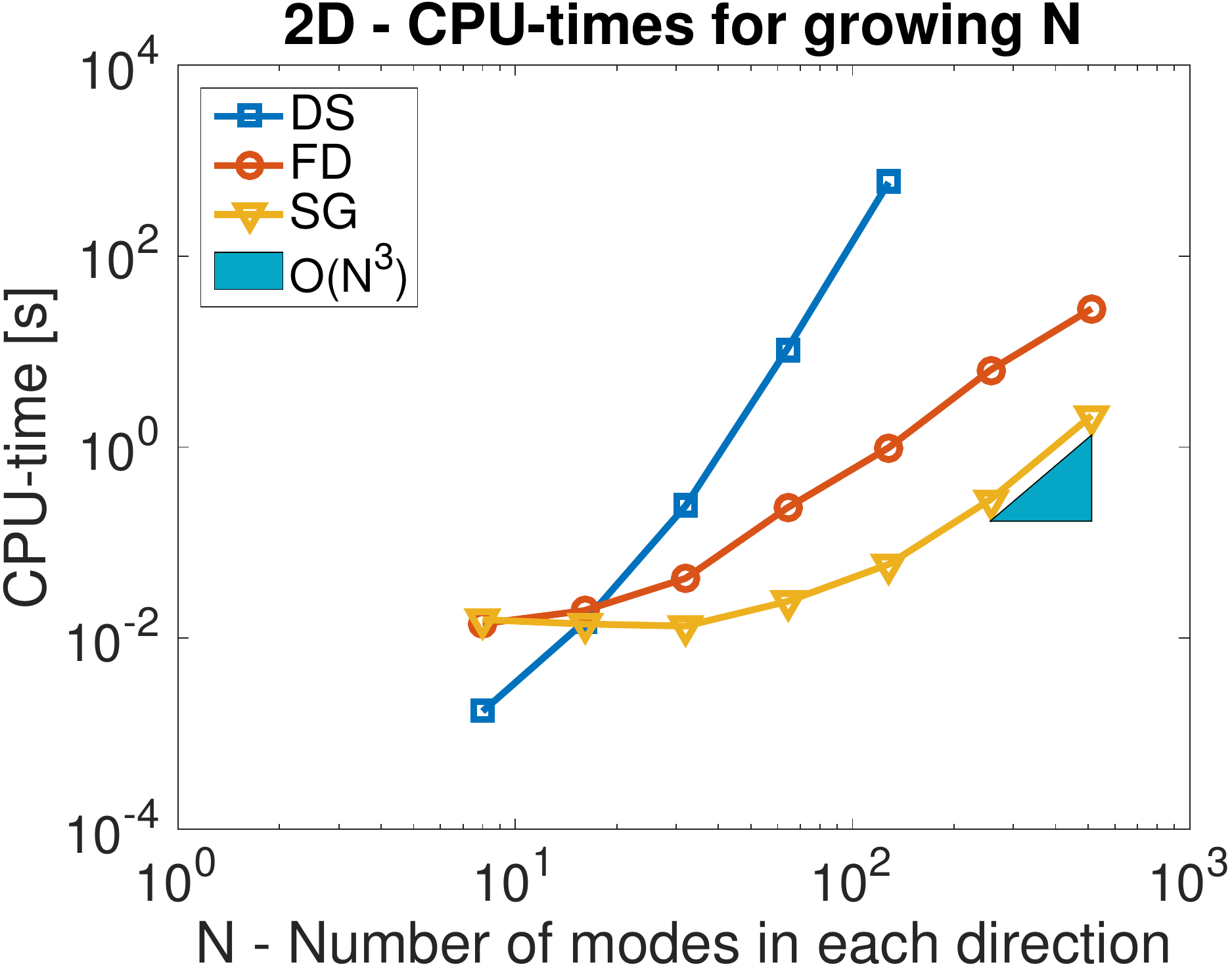}}%
\hfill
\subcaptionbox{3D CPU-times for increasing problem sizes.}{\includegraphics[width=0.45\textwidth]{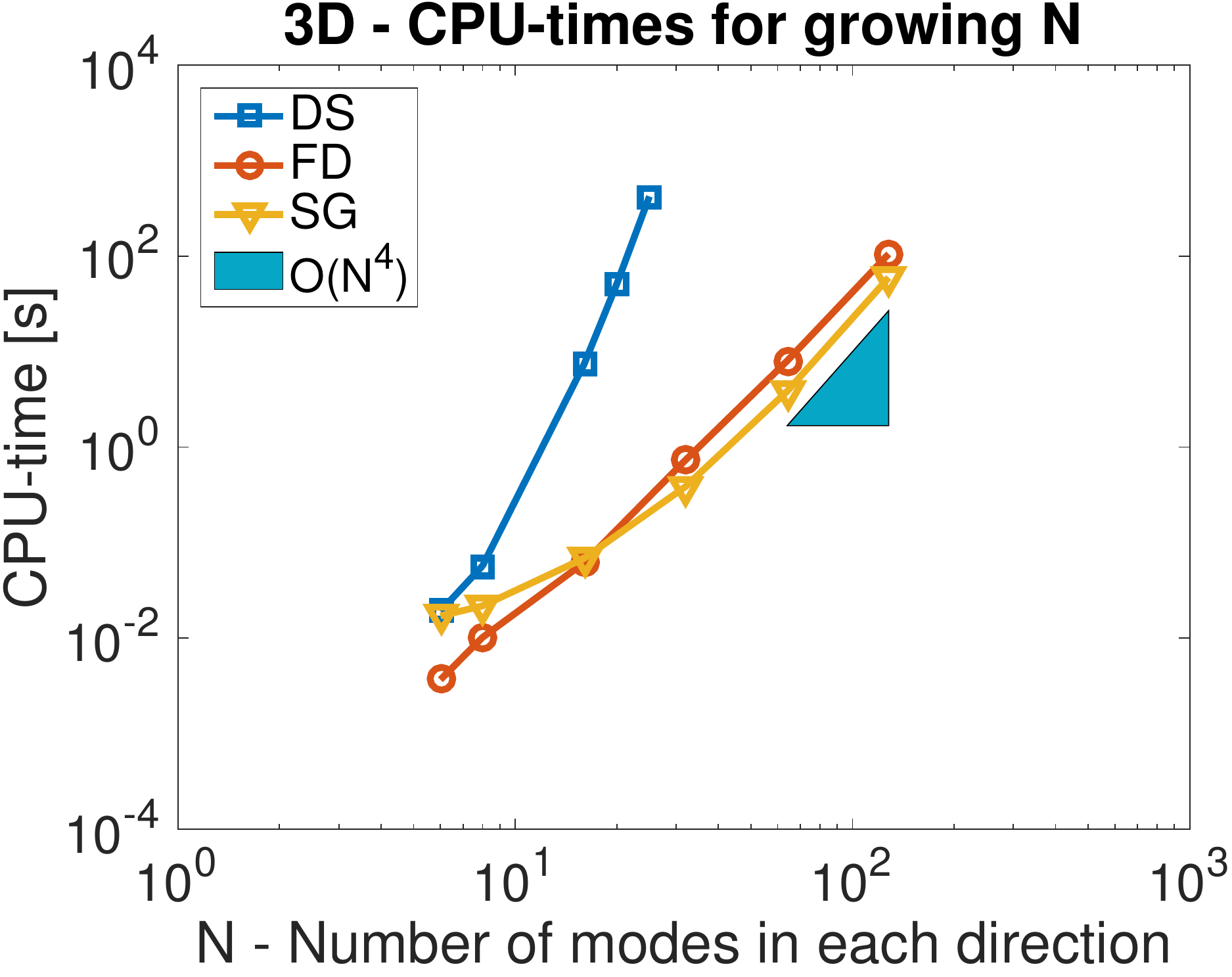}}%
\caption{CPU-times for increasing problem sizes.}
 \label{fig:CPUtime2}
\end{figure}

Case study I considers the CDR model problem (\ref{OCPWeak}) in $\Omega:=[-1,1]^{d}, \hspace{2pt} d=2,3$, where the two-and three dimensional coefficients are defined by 
\begin{equation} a(\bold{x}):=1, \quad \beta=\bold{0}, \quad \gamma(\bold{x})=10^k+\bold{x}, \quad k:=0.8. \end{equation}
Using the continuous optimality conditions (\ref{dualsys}), the desired states, $z_{d} \in L^2(\Omega),$ and the source-terms, $f \in L^2(\Omega),$ have been constructed such that the optimal solutions, $(y^{\ast}, u^{\ast}) \in H_{0}^{1}(\Omega) \times L^2(\Omega),$ are given by 
\begin{alignat}{3} 
   &2\text{D}:  \quad y^{\ast}(\bold{x}) &&=\sin(\pi x_1) \sin(\pi x_2), \quad && u^{\ast}(\bold{x}) = y^{\ast}(\bold{x}) (2 \pi^2+\gamma(\bold{x})), \\
   &3\text{D}:  \quad y^{\ast}(\bold{x}) &&= \sin(\pi x_1) \sin(\pi x_2) \sin(\pi x_3) , \quad && u^{\ast}(\bold{x}) =   y^{\ast}(\bold{x})  (3 \pi^2+\gamma(\bold{x})).
\end{alignat}
The study compares the SG scheme to 1) \texttt{MATLAB}s direct solver (DS) and 2) the second-order FD scheme in 2D and 3D, respectively. Figure \ref{fig:Convergence} illustrates convergence rates, table \ref{tab:2D3D_epsilon} and figure \ref{fig:CPUtime1} compare the CPU-times that each method requires to reach a given level of accuracy while table \ref{tab:2D_CPU} and figure \ref{fig:CPUtime2} show CPU-times for increasing problem sizes. The results confirm exponential convergence rates for the SG scheme. In particular, the SG scheme reaches an accuracy of $10^{-4}$ using 10 modes in each direction. In contrast, the 3D FD scheme requires 128 nodes in each direction to reduce the error to about $4.8 \cdot10^{-4}.$ As table \ref{tab:2D3D_epsilon} shows, this difference in accuracy manifests itself in the CPU-times required to resolve the problem to a given precision. While the SG scheme manages to solve the 3D problem to an error tolerance of $10^{-9}$ in approximately $0.07$ seconds, the FD scheme becomes computationally intractable. Note also that the direct method is slightly faster than the SG scheme in the 2D case. This is due to the low dimensionality that comes with spectral accuracy. However, the direct method becomes inferior in 3D. In particular, it requires approximately $6.4$ seconds of CPU-time to reach an error tolerance of $10^{-9}.$ \par For this case study, spectral accuracy keeps the problem size relatively low, however, potential applications that involve large systems of coupled CDR equations have significantly more degrees of freedom. To demonstrate the SG schemes potential to handle such applications, table \ref{tab:2D_CPU} and figure \ref{fig:CPUtime2} investigate the development of CPU-times for an increasing number of unknowns. As the CPU-timings show, the SG scheme leads to fast solution of large-scale problems on standard hardware. In the 2D case, problems with approximately $130000$ unknowns are solved in less than a half a second. In the 3D case, problems with up to approximately 4 million unknowns are solved in less than a minute. For problems of these sizes, \texttt{MATLAB}s direct solver runs out of memory. Furthermore, the direct method quickly becomes impractical in terms of CPU-time. Overall, the results demonstrate that the SG scheme can produces high-precision solutions to the optimal control problem (\ref{OCPWeak}), while its computational complexity is comparable to that of a simple FD implementation, when a Schur-complement preconditioner is used. 

\begin{table}[tb!]
\centering
\caption{\small 2D and 3D cases. Comparison of CPU-times to reach a given error tolerance, $\epsilon$. An asterix, $\ast$, indicates that the desired level of accuracy was attained by the listing in the above slot. A vertical line, $-$, indicates that no computations were performed.}
\begin{tabular}{| l*{1}{c}|r  l*{2}{c}|r  l*{2}{c} |} 
\hline
& $\epsilon$ & 2D & DS & FD & SG & 3D & DS & FD & SG \\
\hline
& $10^{-1}$ & & $0.0013$ & $0.0391$  & $0.0063 $ &&  $0.0206$ & $0.0752 $  & $0.0207 $ \\
& $10^{-2}$ & & $0.0014$ & $0.4452 $  & $0.0066 $ &&  $0.1066$ & $0.4204 $ & $0.0353 $ \\
& $10^{-3}$ & & $\ast$ & $5.9817 $  & $\ast $ &&  $0.2200$ & $17.791 $  & $0.0373 $ \\
& $10^{-4}$ & & $0.0019$ & $69.172 $  & $0.0066 $ &&  $\ast$ & $1195.3$ & $\ast$ \\
& $10^{-5}$ & & $0.0036$ & $-$  & $0.0069 $ &&  $0.7448$ & $-$ & $0.0438 $ \\
& $10^{-6}$ & & $\ast$ & $-$  & $\ast$ &&  $\ast$ & $-$ & $\ast $\\
& $10^{-7}$ & & $0.0066$ & $-$  & $0.0074$ &&  $2.2466$ & $-$ & $0.0579$\\
& $10^{-8}$ & & $\ast$ & $-$  & $\ast $ &&  $\ast$ & $-$ & $\ast $\\
& $10^{-9}$ & & $0.0114$ & $-$  & $0.0079 $ &&  $6.3949$ & $-$ & $0.0721$\\
\hline 
\end{tabular} \label{tab:2D3D_epsilon} 
\end{table}

\begin{table}[tb!]
\centering
\caption{\small 2D and 3D cases. Comparison of iteration counts and CPU-times for increasing $N$. Iteration counts are listed in parentheses. A vertical line, $-$, indicates that no computations were performed.}
\begin{tabular}{| l*{5}{c} |} 
\hline
& 2D &  &  &  &\\
\hline
& N & DoF & DS & FD & SG \\
\hline
& 8&  $98$ & $0.0017$ & $0.0156 \quad (5)$  & $0.0126 \quad (7)$\\
& 16 & $450$ & $0.0147$ & $0.0194 \quad (5)$ & $0.0141 \quad (7)$ \\
& 32 & $1922$  & $0.0243$ & $0.0423 \quad (5)$  & $0.0133 \quad (7)$\\
& 64 & $7938$  & $10.581$ & $0.2312 \quad (5)$ & $0.0242 \quad (7)$ \\
& 128 & $32258$  & $590.52$ & $0.9878 \quad (5)$ & $0.0593 \quad (7)$ \\
& 256 & $130050$  & $-$ & $6.3912 \quad (5)$ & $0.2821 \quad (7)$\\
& 512 & $522242 $  & $-$ & $23.741 \quad (5)$ & $2.1451 \quad (7)$\\
\hline 
& 3D &  &  &  & \\ 
\hline
& N & DoF & DS & FD & SG \\
\hline
& 6&  $250$ & $0.02$ & $0.004 \quad (2)$  & $0.017 \quad (8)$\\
& 8 & $686 $ & $0.06$ & $0.012 \quad (2)$ & $0.022 \quad (7)$ \\
& 16 & $6750$  & $ 7.37$ & $0.064 \quad (2)$  & $0.068 \quad (7)$\\
& 25 & $27648$  & $ 420.9$ & $0.341 \quad (2)$ & $0.191 \quad (7)$ \\
& 32 & $59582$  & $ -$ & $0.731 \quad (2)$ & $0.383 \quad (7)$ \\
& 64 & $500094$  & $-$ & $8.124 \quad (2)$ & $3.842 \quad (7)$\\
& 128 & $4096766 $  & $-$ & $105.2 \quad (3)$ & $59.58 \quad (7)$\\
\hline 
\end{tabular} \label{tab:2D_CPU} 
\end{table}

\begin{figure}
\centering
\subcaptionbox{Case I. Sensitivity to growing problem sizes and the Tikhonov parameter, $\rho.$}{\includegraphics[width=0.45\textwidth]{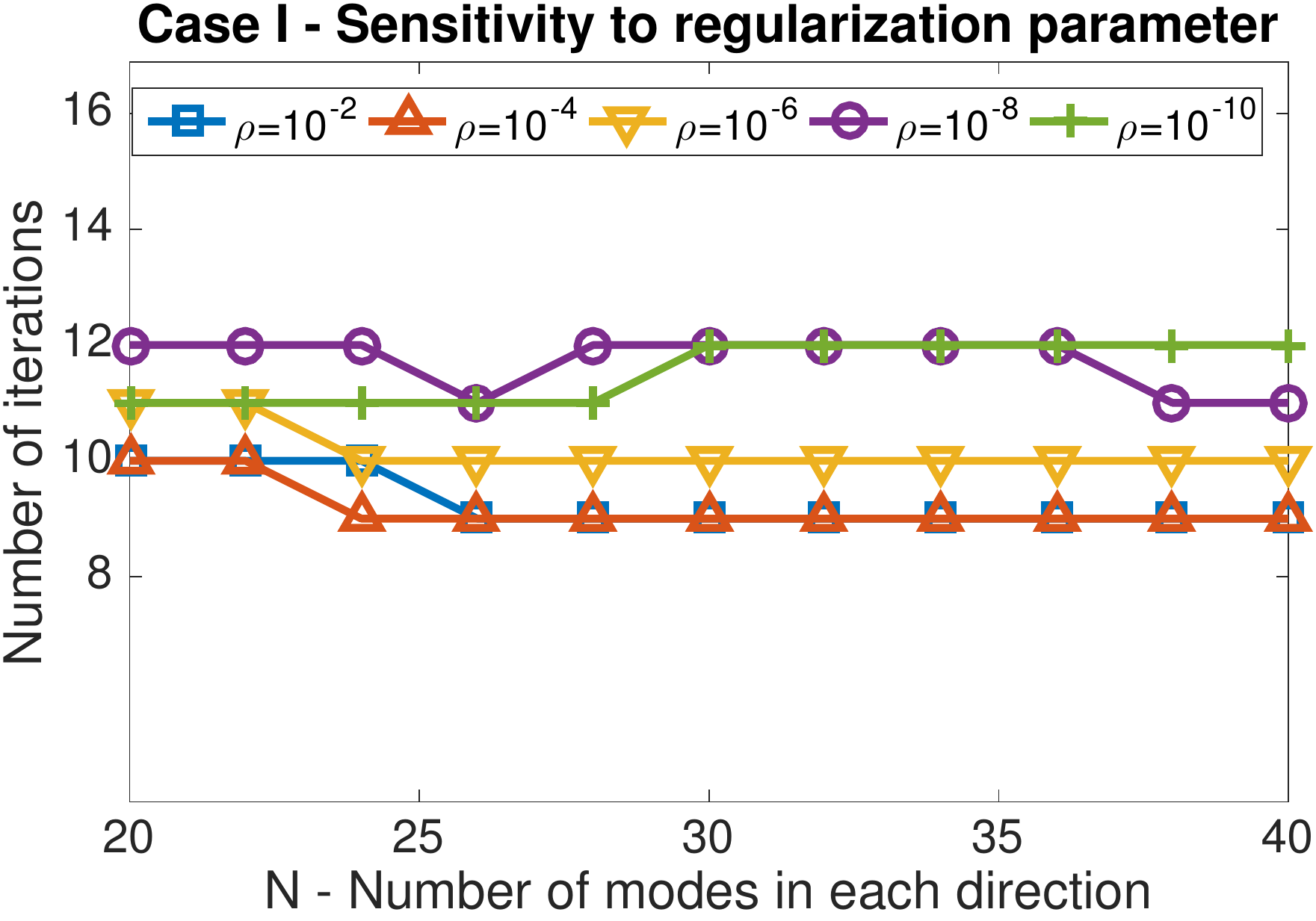}}%
\hfill
\subcaptionbox{Case II. Sensitivity to growing problem sizes and the Tikhonov parameter, $\rho.$}{\includegraphics[width=0.45\textwidth]{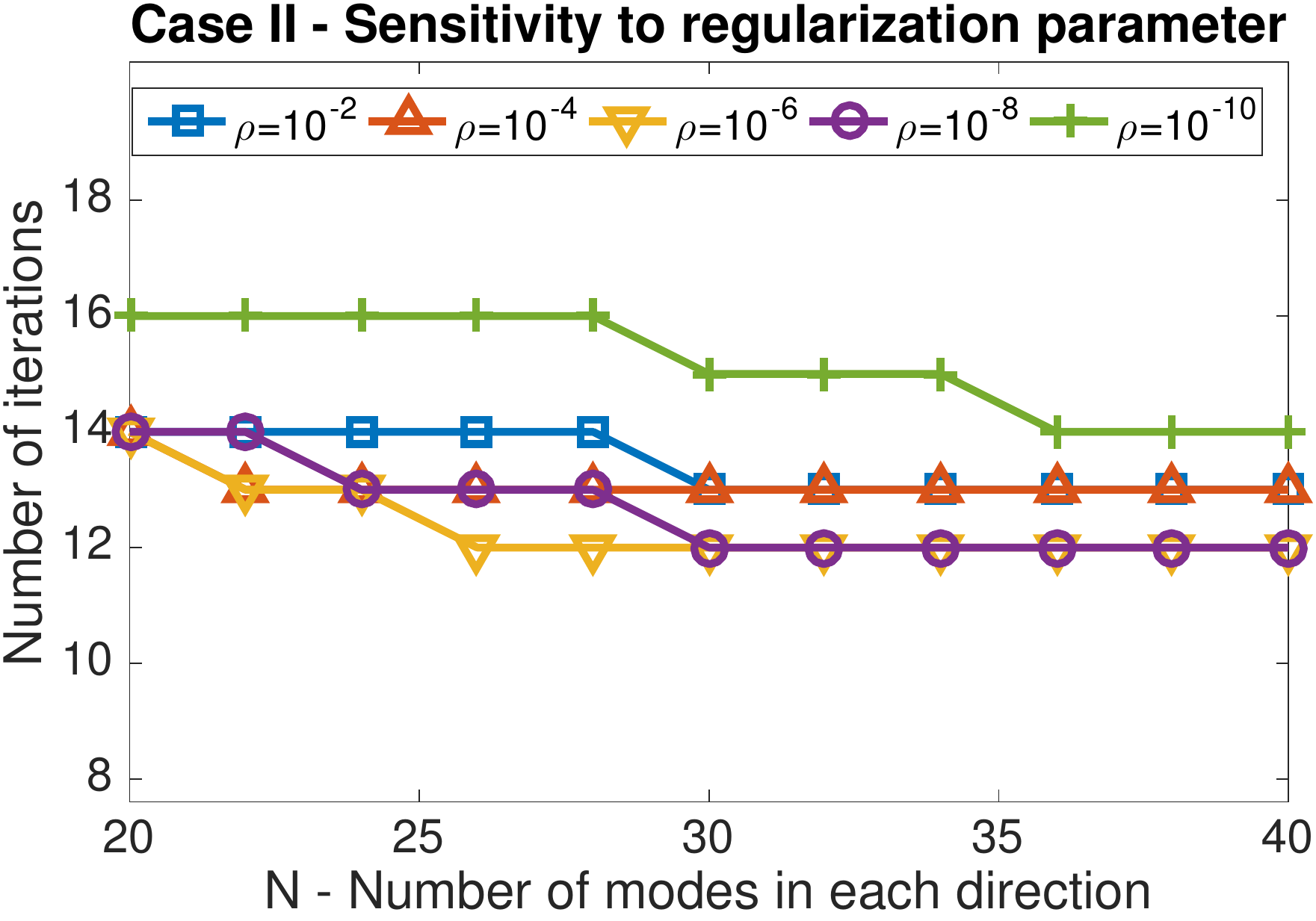}}%

\subcaptionbox{Case I. Eigenvalue distribution of the preconditioned system matrix.}{\includegraphics[width=0.45\textwidth]{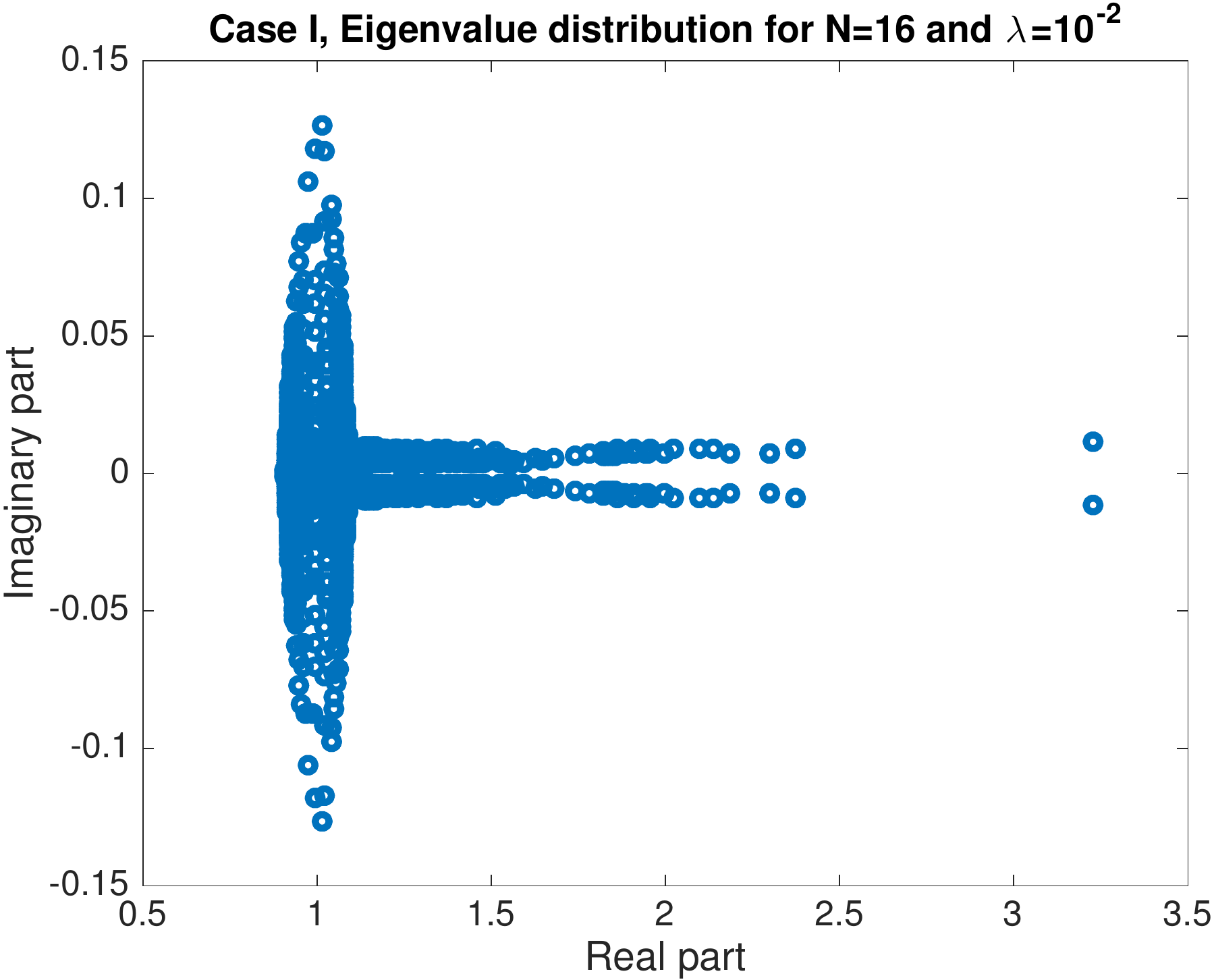}}%
\hfill
\subcaptionbox{Case II. Eigenvalue distribution of the preconditioned system matrix.}{\includegraphics[width=0.45\textwidth]{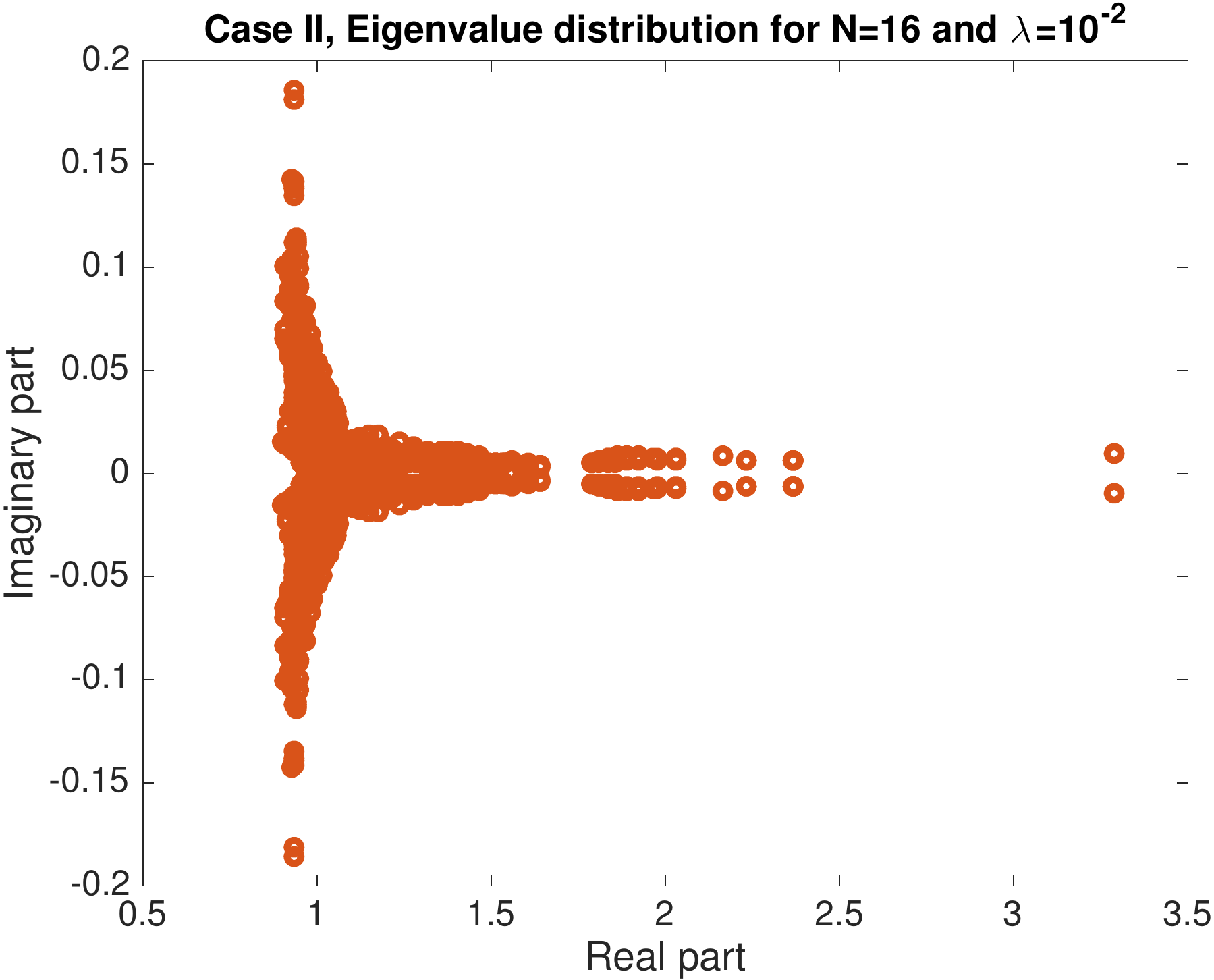}}%
\caption{Iteration counts for increasing problem sizes and the different choices of the Tikhonov parameter, $\rho,$ and eigenvalue distributions.}
\label{fig:iterlambda}
\end{figure}
\subsection{Case study II - Robustness of the SG scheme} 
Case study II treats the 3D CDR model problem (\ref{OCPWeak}) in $\Omega:=[-1,1]^{3}$ for two sets of variable coefficients, $C_1$ and $C_2,$ that are given by:
\begin{alignat}{4} &C_1: \quad a(\bold{x})&&:=10+x_1, \quad &&\beta(\bold{x})=(10+x_1,0,0), \quad &&\gamma(\bold{x})=10+\bold{x}. \\
                     &C_2: \quad a(\bold{x})&&:=10+\bold{x}^2, \quad &&\beta(\bold{x})=(10+x_1^2,10+x_2^2,0), \quad &&\gamma(\bold{x})=10+\bold{x}. 
\end{alignat}

For each of the cases, the corresponding desired states, $z_{d} \in L^2(\Omega),$ and the source-terms, $f \in L^2(\Omega),$ have been chosen such that the optimal solutions, $(y^{\ast}, u^{\ast}) \in H_{0}^{1}(\Omega) \times L^2(\Omega),$ are given by 
\begin{equation}
   y^{\ast}(\bold{x}) = \sin(\pi x_1) \sin(\pi x_2) \sin(\pi x_3), \quad u^{\ast}(\bold{x}) = y^{\ast}(\bold{x}) (3 \pi^2 a(\bold{x})-(\nabla \cdot \beta)(\bold{x})+\gamma(\bold{x})).
\end{equation}
To investigate how the performance of the SG scheme depends on 1) the number of unknowns and 2) the choice of regularization parameter, figure \ref{fig:iterlambda} shows the iteration count for growing $N$ and various choices of $\rho.$ The results show that the iteration counts remain approximately constant as $N$ grows, independently of $\rho.$ This indicates that the preconditioner, $P_{\rho},$ is insensitive to both the problem size and the Tikhonov parameter. Overall, this supports that $P_{\rho}$ is robust for optimal control problems of the type (\ref{OCPWeak}), whenever the coefficients are given by low-order polynomials. 
\par
To provide a heuristic explanation for the effectiveness of the preconditioners, the following considers the eigenvalue distribution associated with the preconditioned system matrices for each of the cases, $C_1$ and $C_2$. In particular, while it is not sufficient to ensure rapid convergence of non-symmetric KSP solvers, tightly clustered eigenvalue distributions are often associated with good performance \cite{saad2003iterative}. For each of the cases, $C_1$ and $C_2$, the subfigures $(\ref{fig:iterlambda} c)$ and $(\ref{fig:iterlambda} d)$ show the respective eigenvalue distributions. In both cases, the eigenvalues are clustered relatively tightly around one and away from zero. Further, experiments that have not been included in the paper, indicate that the clusterings are insensitive to the choice of regularization parameter, $\rho$. Overall, this provides empirical evidence for effectiveness of the preconditioners, $P_{\rho}$.

\section{Conclusions} \label{Sec6}
This paper has proposed a fast and memory-efficient high-order spectral Galerkin (SG) scheme tailored for distributed elliptic optimal control problems governed by convection-diffusion-reaction equations. By design, the scheme targets optimal control problems where the optimization problem is constrained only by the PDE. For this class of problems, the SG scheme converges spectrally fast to the optimal solution provided only that the data of the problem is smooth. Consequently, relative to low-order methods, the SG scheme can significantly reduce computational effort by lowering both storage requirements and CPU-time. At its core, the SG scheme combines a non-symmetric Krylov subspace (KSP) solver with a new block preconditioner tailored for the saddle-point optimality system that arises from the SG discretization. The preconditioner is matrix-free and its complexity scales linearly with the number of modal basis functions, $N^{d}$, where the proportionality constant is small. Overall, the computational complexity of the SG scheme amounts to a small multiple of $N^{d+1}.$ Combined with exponential converges rates and its matrix-free nature, this implies that the SG scheme is competitive, even for most large-scale problems. The computational bottleneck of the SG scheme is tied to the matrix-vector products required by the KSP solver. This observation stresses the importance of an efficient preconditioner that keeps the number of KSP iterations low. The preconditioner proposed by this paper appears to fulfill these requirements for the large class of problems with variable coefficients given by low-order polynomials. In particular, the numerical results strongly indicate that the preconditioner is ideal in the sense that GMRES converged in a low number of iterations, independently of the number of modes and the choice of Tikhonov regularization parameter. In particular, for all cases considered by this paper, the KSP solver converged in less than $16$ iterations. 
Until now, efficient spectral Galerkin schemes only existed for problems with constant coefficients. The SG scheme extends current state-of-the-art spectral Galerkin methods to variable-coefficient problems. Since Newton-like solvers often rely on repeated solution of linearized variable-coefficient problems, the SG scheme can be considered as an important first step towards efficient solution of non-linear optimal control problems. In this way, the SG scheme may broaden applications of high-order methods to PDE-constrained optimal control. 
As a limitation, this paper focuses on problems that are constrained only by the PDE. However, many practical applications also impose integral and bound constraints on the state and control variables. Future work intends to use the SG scheme as an important building block to handle more complex optimal control problems.

\bibliographystyle{siamplain}
\bibliography{References_new}
\end{document}